\documentclass[12pt]{article}
\makeatletter \oddsidemargin  -.1in \evensidemargin -.1in
\newcommand{\singlespacing}{\let\CS=\@currsize\renewcommand{\baselinestretch}{1}\tiny\CS}
\newcommand{\oneandahalfspacing}{\let\CS=\@currsize\renewcommand{\baselinestretch}{1.25}\tiny\CS}
\newcommand{\doublespacing}{\let\CS=\@currsize\renewcommand{\baselinestretch}{1.35}\tiny\CS}

\newtheorem{rule-def}[theorem]{Rule}

\RequirePackage[dvips]{graphicx} \textheight 22.5cm
\usepackage{latexsym,epsfig,enumerate,amsmath,amsfonts,amssymb,amsbsy,amsopn,mathrsfs}
\usepackage{subfigure}
\usepackage{graphicx}

\begin{document}

\title{\bf Multiple and Complete New Important Conjectures on Perfect Cuboid and Euler Brick} \author{\small
  Somnath Maiti$^{1,2,3,}$\thanks{Email address of Corresponding author: {\it
      maiti0000000somnath@gmail.com
      (Somnath Maiti)}}\\\it $^{1}$ Department of Mathematics, Rajendra College \\\it (a constituent unit of Jai Prakash
University), Chapra, Bihar, India\\ \it
  $^{2}$Department of Mathematical Sciences, Indian Institute of
  Technology (BHU),\\\it Varanasi-221005, India\\ \it
  $^{3}$Current Address:- Somnath Maiti, Flat No. A-3, Balaji Vihar, \\\it Village: Rupa Ki Nangal, Post Office: Sumel (Jamdoli),
\\\it District: Jaipur, State: Rajasthan, India, Pin: 302031.}
\date{Research work year: 2023-2024}
\maketitle \noindent \doublespacing
\vspace{-0.5cm}
\begin{abstract}
  Nobody has discovered any perfect cuboid and there is no formula to
  deliver all possible Euler bricks. During investigations of famous
  open problems regarding the perfect cuboid and Euler brick; I have
  found new important conjectures on Pythagorean triples and
  biquadratic Diophantine equations \cite{Maiti2024} which are reduced
  $\&$ complete form for perfect cuboid and Euler brick problems. The
  details of the conjectures have been provided in Sections
  \ref{pythagorean_conjecture}-\ref{biquadratic_conjecture}. If any
  perfect cuboid exists, it will be only among the solutions of six
  conjectures and all the Euler bricks are only among the solutions of
  next three conjectures \cite{Maiti2024}. For example, if any odd
  $n\in \mathbb{N}$ satisfy $n=e^2-f^2=g^2-h^2=k^2-l^2$ and
  $e^2f^2=g^2h^2+k^2l^2$; then we can discover a perfect cuboid of
  type 1 as $\{e^2-f^2,2gh,2kl,g^2+h^2,k^2+l^2,2ef,e^2+f^2\}$ having
  $(e^2-f^2,2gh,2kl)$ as its edges; $(g^2+h^2,k^2+l^2,2ef)$ as its
  face diagonals and $e^2+f^2$ as its body diagonal where
  $e,f,g,h,k,l~(>1)\in \mathbb{N}$. Equivalently, biquadratic
  Diophantine equation conjectures have been introduced for these
  perfect cuboid conjectures. For the benefit of readers, along with
  the original contribution for new important conjectures on perfect
  cuboid and Euler brick problems; brief review related to Pythagorean
  Triple, perfect cuboid and Euler brick problems as well as on
  Diophantine Equation and Biquadratic Diophantine Equation; studied
  in the past by previous researchers, have been discussed in the
  paper.

  \vspace{0.5cm} \it Keywords: {\small  Perfect
    Cuboid Conjecture; Pythagorean Triple Conjecture;
    Euler Brick Conjecture; Biquadratic Diophantine Equation Conjecture; Number Theory Conjecture.}
\end{abstract}

\section{Introduction}
Pythagorean triples \cite{Agarwal2020} $\{x,y,z\}$ are natural numbers
$x,y,z>1$ which satisfy the equation $x^2+y^2=z^2$ and an ordered
Pythagorean triple $(x,y,z)$ is a triple that also satisfies the order
relation $x<y<z$. For any natural number $r$ and any Pythagorean
triple $(x,y,z)$; the triple $\{rx,ry,rz\}$ is also a Pythagorean
triple. A Pythagorean triple $\{(x,y,z\}$ is called primitive if
$gcd(x,y,z)=1$. All primitive Pythagorean triples $\{x,y,z\}$ can be
given parametrically as $$x=u^2-v^2,~y=2uv,~z=u^2+v^2;$$ where
$u,v(>0)$ are natural numbers of opposite parity, $u>v$ and
gcd$(u,v)=1$. Thus anyone can easily deduce that all Pythagorean
triples $(x,y,z)$ can be characterized by
$\{x,y,z\}=\{r(u^2-v^2),r(2uv),r(u^2+v^2)\}$, where $u,v(>0)$ are
natural numbers of opposite parity, $u>v$, $gcd(u,v)=1$ and
$r=gcd(x,y,z)$. All these results are well-known and can be found in
most books on elementary Number Theory. In recent years, researchers
have been using primitive Pythagorean triples in cryptography as random
sequences and for the generation of keys \cite{Kak2014}.

The parametric solution of $x^2+y^2=z^2$ can support in suggesting and
proving several properties that Pythagorean triples satisfy
\cite{Agarwal2020}. For example, 60 divides $xyz$ for each Pythagorean
triple $\{x,y,z\}$. Without difficulty, we can show that for every
$n\ge 3$, there is a Pythagorean triple $\{n,y,z\}$ or $\{x,n,z\}$ or
$\{x,y,n\}$. Moreover, for each $n\ge 1$, there are at least $n$
Pythagorean triples having the same least
member \cite{Agarwal2020,Tripathi2008-09}.

An Euler brick, also called a rational cuboid, is a rectangular
parallelepiped with integer sides and integer face diagonals. As
available in the literature, Paul Halke \cite{Dickson1966} obtained
the first and smallest Euler Brick with side lengths $(a,b,c)=(44,
117, 240)$ and face diagonals $(e,f,g)=(125,244,267)$ in 1719. If $a,
b, c$ are the edges and $d, e, f$ are the face diagonals; the Euler
brick can be formulated mathematically as
\begin{eqnarray}
a^2+b^2=d^2,~a,b,d\in \mathbb{N};
\label{Euler brick 1st}\\
b^2+c^2=e^2,~b,c,e\in \mathbb{N};
\label{Euler brick 2nd}\\
a^2+c^2=f^2,~a,c,f\in \mathbb{N};
\label{Euler brick 3rd}
\end{eqnarray}
Nicholas Saunderson \cite{Saunderson} found a parametric solution to
the Euler Brick (known as the Euler cuboid) giving families of Euler
bricks and Euler introduced at least two parametric solutions in 1770
and 1772, although these do not deliver all possible Euler bricks.

A perfect cuboid is defined as a cuboid where the edges, face
diagonals and space diagonal all have integer lengths. Thus, the
following equation
\begin{eqnarray}
a^2+b^2+c^2=g^2;~a,b,c,g\in \mathbb{N}
\label{Perfect Euler brick 1st}
\end{eqnarray}
can be added to the system of Diophantine equations (\ref{Euler brick
  1st})-(\ref{Euler brick 3rd}) defining an Euler brick.

Does such a perfect cuboid exist? Nobody has discovered any perfect
cuboid nor it has been proved that a perfect cuboid does not
exist. Moreover, as reported by Spohn \cite{Spohn1,Spohn2}, Chein
\cite{Chein}, Lagrange \cite{Lagrange} and \cite{Leech}; either the
Euler cuboid or its derived cuboid can not be a perfect cuboid. Euler
announced an example where only the body diagonal failed to be an
integer value (Euler brick). with rapidly increasing degrees, Colman
\cite{Colman} introduced infinitely many two-parameter
parametrizations of rational cuboids whose all the seven lengths
are integers except possibly for one edge, called an edge cuboid, or
face diagonal, called a face cuboid.

As reported by Kraitchik \cite{Kraitchik1}, at least one of the sides
of a rational cuboid has a divisor 4 and another one is divided by
16. Also, the sides of a perfect cuboid have divisors of different
powers of 3 and at least one of the sides is divisible by both the
primes 5 and 11. Horst Bergmann carried the equally elementary
extension of this result and Leech pointed out that the product of all
the sides and diagonals (edge and face) of a perfect cuboid is
divisible by $2^8\times 3^4\times 5^3\times 7\times 11\times 13\times
17\times 19\times 29\times 37$ (cf. \cite{Guy}, Problem
D18). Kraitchik \cite{Kraitchik1} rediscovered the Euler cuboids given
by Saunderson and provided a list of 50 rational cuboids that are not
Euler cuboids by using some ad hoc methods. He further extended his
classical list to 241 cuboids in \cite{Kraitchik2} in which the odd
sides were less than $10^6$ and also reported 18 more in \cite{Kraitchik3}
giving 16 new ones. Lal and Blundon \cite{Lal} showed that for integers
$m,~n,~p$ and $q$; the cuboid having sides
$x=|2mnpq|,~y=|mn(p^2-q^2)|$ and $z=|pq(m^2-n^2)|$ has at least two
face diagonals as integers and is a rational cuboid if and only if
$y^2+z^2=\square$. Using the symmetry, the computer searches by them through
all the quadruples $(m,~n,~p,~q)$ satisfying $1\leq m, n, p, q\leq 70$
and if $y^2+z^2=\square$; provided 130 rational cuboids which are not
perfect cuboids, however, Shanks \cite{Shanks} found some corrigenda
in their paper.

Korec \cite{Korec1} pointed out that if $x,~y,~z$ are the sides of a
perfect cuboid, then we can find natural numbers $a,~b,~c$ all
dividing $x$ and $t=\sqrt{y^2+z^2}$ such that
$$y=\frac{1}{2}\left(\frac{x^2}{a}-a\right),~z=\frac{1}{2}\left(\frac{x^2}{b}-b\right),~t=\frac{1}{2}\left(\frac{x^2}{c}-c\right)$$. With
this consideration, he reported no perfect cuboids when the least side
smaller than 1000 and he \cite{Korec2} extended his result indicating
no perfect cuboid if the least edge is smaller than $10^6$. Korec
\cite{Korec3} also not found any perfect cuboid if the full diagonal
of a perfect cuboid is $<8\times 10^9$, however, if $x$ and $z$ are
the maximal edge and the full diagonal respectively of a perfect
cuboid, then $z\leq x\sqrt{3}$. If a side $x\leq 333750000$, Rathbun
\cite{Rathbun1} announced 6800 body, 6749 faces, 6380 edges and no
perfect cuboids by a computer search. He \cite{Rathbun2} also pointed
out that 4839 of the 6800 body or rational cuboids have an odd side
less than 333750000 which is an extension and correction of Kraitchik’s
classical table \cite{Kraitchik1,Kraitchik2,Kraitchik3}. Butler
\cite{Butler} declared no perfect cuboids despite a thorough computer
search when all ``odd sides''$\leq 10^{10}$.

Diophantine equations are among the oldest subjects in mathematics
dealing with equations with integer coefficients in which integer
solutions are studied. For any natural number $n$, Leonard Euler
\cite{Euler1772,Dickson1966} conjectured in 1772 that we would always
find an example where the sum of $n$ $n$th powers would be an $n$th
power, but we could not find any less than $n$ $n$th powers sum equal
to be an $n$th power. Hence, he asserted that $a^4+b^4+c^4=d^4$, known
as the Euler Quartic Conjecture, and $a^5+b^5+c^5+d^5=e^5$ had no
solutions in the integers, but $a^4+b^4+c^4+d^4=e^4$ would have an integer
solutions although he did not find any.

A computer search, nearly two centuries later by Lander and
Parkin \cite{Lander1966}, discovered the first counterexample to the
general conjecture for $n=5$ as $27^5+84^5+110^5+133^5=144^5$. The Euler's
conjecture was also disproved by Noam Elkies \cite{Elkies1987-88} in 1987
for $n=4$ by reporting that $a^4+b^4+c^4=d^4$ has infinitely many
solutions using a geometric construction, however, it is not known if
a parametric solution exists \cite{Ahmadi2023}.

If a sum of $m$ quartics equals to a sum of $n$ fourth powers, the
corresponding quartic (biquadratic) Diophantine equation is called
quartic $m-n$ equation. The $2-1$ equation $a^4+b^4=c^4$, a special
case of Fermat’s Last Theorem, has no solution and in general the
equations $a^4\pm b^4=c^2$ have no integer solutions. It is known that
the quartic equation $x^4+py^4=z^4$ has no nonzero integer solutions
\cite{Manley2006}. As reported by several number theorists
(\cite{Choudhry1995,Izadi2017}; \cite{Dickson1992}, pp 648), the
quartic equation $A^4+hB^4=C^4+hD^4$ has a solution in integers for an
infinite number of positive integer values of $h$, however, it is not
known whether the said equation has a solution in integers for all
positive integer values of $h$. Janfada and Nabardi \cite{Janfada2019}
obtained a necessary condition for $n$ to have an integral solution
for the equation $x^4+y^4=n(u^4+v^4)$ in 2019 and reported a
parametric solution. Parametric solutions to the $2-3$ quartic
equation $a^4+b^4+c^4=d^4+e^4$ are known \cite{Ferari1913} and the
smallest solution is $7^4+7^4=3^4+5^4+8^4$. Ramanujan
\cite{Ramanujan1987,Ahmadi2023} discovered the general expression to
the $3-3$ quartic equation $a^4+b^4+c^4=d^4+e^4+f^4$ which can be
represented by
$3^4+(2x^4-1)^4+(4x^5+x)^4=(4x^4+1)^4+(6x^4-3)^4+(4x^5-5x)^4$.

Several people had found solutions of $a^4+b^4+c^4+d^4=e^4$
together with Elkies' infinite family with $d=0$ and Euler essentially
conjectured that this equation has infinitely many primitive integral
solutions. Rose and Brudno listed the first 82 consecutive primitive
solutions of this equation in a table \cite{Rose1973}. Jacobi and
Madden \cite{Jacobi2008} reported that there are infinitely many
primitive integral solutions of this equation with all terms nonzero,
however, the existence of any parametric solution is not known
\cite{Ahmadi2023}. Alvarado and Delorme \cite{Alvarado2014}
established two solutions of $a^4+b^4+c^4+d^4=e^2$ in which $a,b,c,d$
are given by polynomials of degrees 16 and 28 in two parameters. They
also depicted how infinitely many such solutions may be obtained by
using elliptic curves. Choudhry et al. \cite{Choudhry2021} obtained
infinitely many solutions of this equation in terms of quadratic
polynomials in two parameters of $a,b,c,d$. As reported by Choudhry et
al. \cite{Choudhry2020} in 2020, the equation
$(x_1^4+x_2^4)(y_1^4+y_2^4)=z_1^4+z_2^4$ has infinitely many
parametric solutions.

The main purpose of this article is to report multiple new and
important conjectures on Pythagorean triple and biquadratic
Diophantine equation, which are complete $\&$ reduced  to perfect
cuboid and Euler brick problem, to Mathematics as well as
Mathematics-related community.

\section{Pythagorean Triple Conjectures for perfect cuboids and Euler bricks}
\label{pythagorean_conjecture}
Any odd $n\in \mathbb{N}$ can be written as
$n=p_1^{\alpha_1}p_2^{\alpha_1}\cdots
p_m^{\alpha_m}=e_1^2-f_1^2=e_2^2-f_2^2=\cdots
=e_{2^{m-1}}^2-f_{2^{m-1}}^2$ for only $2^{m-1}$ distinct primitive
Pythagorean triples with an odd edge $n$, since for each
$n=p_1^{\alpha_1}p_2^{\alpha_2}\cdots p_m^{\alpha_m}=e_i^2-f_i^2$, we
have a primitive Pythagorean triple
$\{e_i^2-f_i^2,2e_if_i,e_i^2+f_i^2\};~i=1,2,\cdots,2^{m-1}$
\cite{Tripathi2008-09,Maiti2024} since $(e_i^2-f_i^2)^2+4e_i^2f_i^2=(e_i^2+f_i^2)^2$ for all $i$.

All the non-primitive Pythagorean
triples can be generated from the following expression on $n$. Let $z$
(except $1$ and $n$) is a divisor of $n$ with
$\frac{n}{t}=z=p_{1}^{\alpha_{a1}}p_{2}^{\alpha_{a2}}\cdots
p_m^{\alpha_{am}};~0\leq\alpha_{ai}\leq\alpha_{i},~i=1,2,\cdots,m$;
then each primitive Pythagorean triple
$\{e_{ai}^2-f_{ai}^2,2e_{ai}f_{ai},e_{ai}^2+f_{ai}^2\}$ with edge $z$
gives rise to a non-primitive Pythagorean triple
$\{t(e_{ai}^2-f_{ai}^2),2te_{ai}f_{ai},t(e_{ai}^2+f_{ai}^2)\}$ with an
odd edge $z$.  If $z$ has $m_1~(\leq m)$ distinct prime factors, then
we can get $2^{m_ 1-1}$ non-primitive Pythagorean triples with an odd
edge $n$ for the corresponding divisor $z$. Thus, total number of
Pythagorean triples (primitive as well as non-primitive) with an odd
edge ($z=n$ corresponds all primitive Pythagorean triples) $n$
is $$\sum_{\substack{z|n \\(z\neq 1)}}2^{m_ 1-1}.$$ Hence, if a
perfect cuboid exists it will be only among the solutions of the
following six conjectures and all the primitive Euler bricks are only
among the solutions of next three conjectures \cite{Maiti2024}.

\subsection{Pythagorean Conjecture 1: perfect cuboid first type}
\label{pythagorean_conjecture1}
Can you provide any odd $n\in \mathbb{N}$ satisfying
\begin{eqnarray}
  n=e^2-f^2=g^2-h^2=k^2-l^2
  \label{pythagorean_conjecture1.1}\\
\text{and}~~~~~~~~~e^2f^2=g^2h^2+k^2l^2;
\label{pythagorean_conjecture1.2}
\end{eqnarray}
where $e, f, g, h, k, l~(>1)\in \mathbb{N}$? We can note that the
equation (\ref{pythagorean_conjecture1.1}) is always true if $n$ has
at least three distinct prime factors. Such type of one $n$ satisfying
equations
(\ref{pythagorean_conjecture1.1}-\ref{pythagorean_conjecture1.2}) can
produce one primitive perfect cuboid of the first type as
$\{e^2-f^2,2gh,2kl,g^2+h^2,k^2+l^2,2ef,e^2+f^2\}$ having
$(e^2-f^2,2gh,2kl)$ as its edges; $(g^2+h^2,k^2+l^2,2ef)$ as its face
diagonals and $e^2+f^2$ as its body diagonal \cite{Maiti2024}. Here, one can easily verify
that $(e^2-f^2)^2+(2gh)^2=(g^2+h^2)^2$,
$(e^2-f^2)^2+(2kl)^2=(k^2+l^2)^2$, $4g^2h^2+4k^2l^2=4e^2f^2$,
$(e^2-f^2)^2+(2gh)^2+(2kl)^2=(e^2+f^2)^2$.

\subsection{Pythagorean Conjecture 2: perfect cuboid second type}
\label{pythagorean_conjecture2}
Can you discover any odd $n\in \mathbb{N}$ satisfying
\begin{eqnarray}
  n=a(e^2-f^2)=g^2-h^2=k^2-l^2
  \label{pythagorean_conjecture2.1}\\
\text{and}~~~~~~~~~a^2e^2f^2=g^2h^2+k^2l^2;
\label{pythagorean_conjecture2.2}
\end{eqnarray}
where $a,e, f, g, h, k, l~(>1)\in \mathbb{N}$? The equation
(\ref{pythagorean_conjecture2.1}) is always true if $n$ has at least
two distinct prime factors and one prime factor of $n$ occurs more
than one time. Similarly, every such $n$ satisfying equations
(\ref{pythagorean_conjecture2.1}-\ref{pythagorean_conjecture2.2}) can
generate one perfect cuboid of the second type as
$\{a(e^2-f^2),2gh,2kl,g^2+h^2,k^2+l^2,2aef,a(e^2+f^2)\}$ having
$(a(e^2-f^2),2gh,2kl)$ as its edges; $(g^2+h^2,k^2+l^2,2aef)$ as its
face diagonals and $a(e^2+f^2)$ as its body diagonal
\cite{Maiti2024}. Here, we can satisfy that
$a^2(e^2-f^2)^2+(2gh)^2=(g^2+h^2)^2$,
$a^2(e^2-f^2)^2+(2kl)^2=(k^2+l^2)^2$, $4g^2h^2+4k^2l^2=4a^2e^2f^2$,
$a^2(e^2-f^2)^2+(2gh)^2+(2kl)^2=a^2(e^2+f^2)^2$.

\subsection{Pythagorean Conjecture 3: perfect cuboid third type}
\label{pythagorean_conjecture3}
Can you find any odd $n\in \mathbb{N}$ so that
\begin{eqnarray}
  n=e^2-f^2=b(g^2-h^2)=k^2-l^2
  \label{pythagorean_conjecture3.1}\\
\text{and}~~~~~~~~~e^2f^2=b^2g^2h^2+k^2l^2;
\label{pythagorean_conjecture3.2}
\end{eqnarray}
where $b,e, f, g, h, k, l~(>1)\in \mathbb{N}$? Such kind of $n$
satisfying equations
(\ref{pythagorean_conjecture3.1}-\ref{pythagorean_conjecture3.2}) can
provide one perfect cuboid of the third type as
$\{(e^2-f^2),2bgh,2kl,b(g^2+h^2),k^2+l^2,2ef,(e^2+f^2)\}$ having
$((e^2-f^2),2bgh,2kl)$ as its edges; $(b(g^2+h^2),k^2+l^2,2ef)$ as its face
diagonals and $(e^2+f^2)$ as its body diagonal \cite{Maiti2024} since  $(e^2-f^2)^2+(2bgh)^2=b^2(g^2+h^2)^2$,
$(e^2-f^2)^2+(2kl)^2=(k^2+l^2)^2$, $4b^2g^2h^2+4k^2l^2=4e^2f^2$,
$(e^2-f^2)^2+(2bgh)^2+(2kl)^2=(e^2+f^2)^2$

\subsection{Pythagorean Conjecture 4: perfect cuboid fourth type}
\label{pythagorean_conjecture4}
Is there any odd $n\in \mathbb{N}$ satisfying
\begin{eqnarray}
  n=e^2-f^2=b(g^2-h^2)=c(k^2-l^2)
  \label{pythagorean_conjecture4.1}\\
\text{and}~~~~~~~~~e^2f^2=b^2g^2h^2+c^2k^2l^2;
\label{pythagorean_conjecture4.2}
\end{eqnarray}
where $b,c,e, f, g, h, k, l~(>1)\in \mathbb{N}$? The equation
(\ref{pythagorean_conjecture4.1}) is always true if either at least
one prime factor of $n$ occurs multiple times or at least two prime
factors of $n$ occur more than one time. Every Such kind of $n$
satisfying equations
(\ref{pythagorean_conjecture4.1}-\ref{pythagorean_conjecture4.2}) can
provide one perfect cuboid of the fourth type as
$\{(e^2-f^2),2bgh,2ckl,b(g^2+h^2),c(k^2+l^2),2ef,(e^2+f^2)\}$ having
$((e^2-f^2),2bgh,2ckl)$ as its edges; $(b(g^2+h^2),c(k^2+l^2),2ef)$ as
its face diagonals and $(e^2+f^2)$ as its body diagonal
\cite{Maiti2024} satisfying  $(e^2-f^2)^2+(2bgh)^2=b^2(g^2+h^2)^2$,
$(e^2-f^2)^2+(2ckl)^2=c^2(k^2+l^2)^2$, $4b^2g^2h^2+4c^2k^2l^2=4e^2f^2$,
$(e^2-f^2)^2+(2bgh)^2+(2ckl)^2=(e^2+f^2)^2$.

\subsection{Pythagorean Conjecture 5: perfect cuboid fifth type}
\label{pythagorean_conjecture5}
Is there any odd $n\in \mathbb{N}$ such that
\begin{eqnarray}
  n=a(e^2-f^2)=b(g^2-h^2)=k^2-l^2
  \label{pythagorean_conjecture5.1}\\
\text{and}~~~~~~~~~a^2e^2f^2=b^2g^2h^2+k^2l^2;
\label{pythagorean_conjecture5.2}
\end{eqnarray}
where $a,b,e, f, g, h, k, l~(>1)\in \mathbb{N}$? The equation
(\ref{pythagorean_conjecture5.1}) is always true if either at least
one prime factor of $n$ occurs multiple times or at least two prime
factors of $n$ occur more than one time. If any odd $n$ satisfy
equations
(\ref{pythagorean_conjecture5.1}-\ref{pythagorean_conjecture5.2}), it
can produce one perfect cuboid of the fifth type as
$\{a(e^2-f^2),2bgh,2kl,b(g^2+h^2),k^2+l^2,2aef,a(e^2+f^2)\}$ having
$(a(e^2-f^2),2bgh,2kl)$ as its edges; $(b(g^2+h^2),k^2+l^2,2aef)$ as its face
diagonals and $a(e^2+f^2)$ as its body diagonal \cite{Maiti2024} holding $a^2(e^2-f^2)^2+(2bgh)^2=b^2(g^2+h^2)^2$,
$a^2(e^2-f^2)^2+(2kl)^2=(k^2+l^2)^2$, $4b^2g^2h^2+4k^2l^2=4a^2e^2f^2$,
$a^2(e^2-f^2)^2+(2bgh)^2+(2kl)^2=a^2(e^2+f^2)^2$.

\subsection{Pythagorean Conjecture 6: perfect cuboid sixth type}
\label{pythagorean_conjecture6}
Can you explore any odd $n\in \mathbb{N}$ which satisfies
\begin{eqnarray}
  n=a(e^2-f^2)=b(g^2-h^2)=c(k^2-l^2)
  \label{pythagorean_conjecture6.1}\\
\text{and}~~~~~~~~~a^2e^2f^2=b^2g^2h^2+c^2k^2l^2;
\label{pythagorean_conjecture6.2}
\end{eqnarray}
where $a,b,c,e, f, g, h, k, l~(>1)\in \mathbb{N}$? The equation
(\ref{pythagorean_conjecture6.1}) is always true if either at least
one prime factor of $n$ occurs multiple times or at least three prime
factors of $n$ occur more than one time. Thus every such $n$
satisfying equations
(\ref{pythagorean_conjecture6.1}-\ref{pythagorean_conjecture6.2}) can
help us to discover one perfect cuboid of the sixth type as
$\{a(e^2-f^2),2bgh,2ckl,b(g^2+h^2),c(k^2+l^2),2aef,a(e^2+f^2)\}$
having $(a(e^2-f^2),2bgh,2ckl)$ as its edges;
$(b(g^2+h^2),c(k^2+l^2),2aef)$ as its face diagonals and $a(e^2+f^2)$
as its body diagonal \cite{Maiti2024} satisfying the equations
$a^2(e^2-f^2)^2+(2bgh)^2=b^2(g^2+h^2)^2$,
$a^2(e^2-f^2)^2+(2ckl)^2=c^2(k^2+l^2)^2$,
$4b^2g^2h^2+4c^2k^2l^2=4a^2e^2f^2$,
$a^2(e^2-f^2)^2+(2bgh)^2+(2ckl)^2=a^2(e^2+f^2)^2$.

\subsection{Pythagorean Conjecture 7: Euler
bricks first type}
\label{pythagorean_conjecture7}
Find all odd $n\in \mathbb{N}$ satisfying
\begin{eqnarray}
  n=e^2-f^2=g^2-h^2 ~\text{and}~e^2f^2+g^2h^2=d^2;
  \label{pythagorean_conjecture7.1}
\end{eqnarray}
where $d,e, f, g, h, ~(>1)\in \mathbb{N}$. Moreover, you can also find
maximum number of occurrence of $[\{(e,f),(g,h)\},~\{ef,gh,d\}]$ for
any particular $n$ and hence maximum occurrence for any $n$. Every
such $n$, satisfying equation (\ref{pythagorean_conjecture7.1}), can
generate one primitive Euler brick of the first type as
$\{(e^2-f^2),2ef,2gh,(e^2+f^2),(g^2+h^2),2d\}$ having
$((e^2-f^2),2ef,2gh)$ as its edges; $((e^2+f^2),(g^2+h^2),2d)$ as its
face diagonals \cite{Maiti2024}. We can easily verify that
$(e^2-f^2)^2+(2ef)^2=(e^2+f^2)^2$,
$(e^2-f^2)^2+(2gh)^2=(g^2+h^2)^2$, $4e^2f^2+4g^2h^2=4d^2$.

\subsubsection{Example:}
It has been verified that there is no Euler bricks first type with the
smallest odd edge less than 1000 in a Table given by Rathbun
\cite{Rathbun3}.

\subsection{Pythagorean Conjecture 8: Euler
bricks second type}
\label{pythagorean_conjecture8}
Provide all odd $n\in \mathbb{N}$ satisfying
\begin{eqnarray}
  n=b(e^2-f^2)=(g^2-h^2)~\text{and}~e^2f^2+b^2g^2h^2=d^2
  \label{pythagorean_conjecture8.1}
\end{eqnarray}
where $b,d,e, f, g, h~(>1)\in \mathbb{N}$. You can also search maximum
number of occurrence of $$[\{(a,e,f),(g,h)\},~\{aef,gh,d\}]$$ for any
particular $n$ and hence maximum occurrence for any $n$. Similarly,
such kind of every $n$, satisfying equation
(\ref{pythagorean_conjecture7.1}), can produce one primitive Euler
brick of the second type as
$\{(e^2-f^2),2ef,2bgh,(e^2+f^2),b(g^2+h^2),2d\}$ having
$((e^2-f^2),2ef,2bgh)$ as its edges; $((e^2+f^2),b(g^2+h^2),2d)$ as its face
diagonals \cite{Maiti2024}.

\subsubsection{Example:}
(i) We have
$85=5\times
17=43^2-42^2=11^2-6^2=5(9^2-8^2)=17(3^2-2^2)$. But
only $85=11^2-6^2=5(9^2-8^2)$ gives
$11^26^2+5^29^28^2=366^2=d^2$. Thus we have the smallest Euler bricks
second type $(85,132,720,157,725,732)$ for an odd edge, but the
smallest Euler brick for any edge is $(117,44,240,125,267,244)$.

(ii) We can calculate
$117=3^2\times 13=59^2-58^2=11^2-2^2=3(20^2-19^2)=3(8^2-5^2)=9(7^2-6^2)=13(5^2-4^2)=39(2^2-1^2)$. But
only $117=11^2-2^2=3(8^2-5^2)$ gives
$11^22^2+3^28^25^2=122^2=d^2$. Thus we have the 2nd smallest Euler bricks
second type $(117,44,240,125,267,244)$ for odd edge, however, it is the
smallest Euler brick for even edge or any edge.

(iii) We evaluate as
$195=3\times 5\times
13=98^2-97^2=34^2-31^2=22^2-17^2=14^2-1^2=3(33^2-32^2)=3(9^2-4^2)=5(20^2-19^2)=5(8^2-5^2)=13(8^2-7^2)=13(4^2-1^2)=15(7^2-6^2)=39(3^2-2^2)=65(2^2-1^2)$. But
only $195=22^2-17^2=3(33^2-32^2)$ gives
$22^217^2+3^233^232^2=3190^2=d^2$. Then, we have the next smallest Euler
bricks second type $(195,748,6336,773,6339,6380)$ for an odd edge.

(iv) We obtain
$231=3\times 7\times
11=116^2-115^2=40^2-37^2=20^2-13^2=16^2-5^2=3(39^2-38^2)=3(9^2-2^2)=7(17^2-16^2)=7(7^2-4^2)=11(11^2-10^2)=11(5^2-2^2)=21(6^2-5^2)=33(4^2-3^2)=77(2^2-1^2)$. But
only $231=16^2-5^2=33(4^2-3^2)$ gives
$16^25^2+33^24^23^2=404^2=d^2$. Then, we have the next smallest Euler
bricks second type as $(231,160,792,281,825,808)$ for an odd edge.

(v) We obtain
$275=5^2\times 11=138^2-137^2=18^2-7^2=5(28^2-27^2)=5(8^2-3^2)=11(13^2-12^2)=25(6^2-5^2)=55(3^2-2^2)$. But
only $275=18^2-7^2=5(8^2-3^2)$ gives
$18^27^2+5^28^23^2=174^2=d^2$. Then, we have the next smallest Euler
bricks second type as $(275,252,240,373,365,348)$ for an odd edge.

(vi) We obtain $495=3^2\times 5\times
11=248^2-247^2=52^2-47^2=32^2-23^2=28^2-17^2=3(83^2-82^2)=3(29^2-26^2)=3(19^2-14^2)=3(13^2-2^2)=5(50^2-49^2)=5(10^2-1^2)=9(28^2-27^2)=9(8^2-3^2)=11(23^2-22^2)=11(7^2-2^2)=15(17^2-16^2)=15(7^2-4^2)=33(8^2-7^2)=33(4^2-1^2)=45(6^2-5^2)=55(5^2-4^2)=99(3^2-2^2)$. But
only $495=52^2-47^2=15(17^2-16^2)$ gives
$52^247^2+15^217^216^2=4756^2=d^2$. Then, we have the next smallest Euler
bricks second type as $(495,4888,8160,4913,8175,9512)$ for an odd edge.

(vii) We obtain $855=3^2\times 5\times
19=428^2-427^2=90^2-85^2=52^2-43^2=32^2-13^2=3(143^2-142^2)=3(49^2-46^2)=3(31^2-26^2)=3(17^2-2^2)=5(86^2-85^2)=5(14^2-5^2)=9(48^2-47^2)=9(12^2-7^2)=15(29^2-28^2)=15(11^2-8^2)=19(23^2-22^2)=19(7^2-2^2)=45(10^2-9^2)=57(8^2-7^2)=57(4^2-1^2)=95(5^2-4^2)=171(3^2-2^2)$. But
only $855=32^2-13^2=15(11^2-8^2)$ gives
$32^213^2+15^211^28^2=1384^2=d^2$. Then, we have the next smallest Euler
bricks second type as $(855,832,2640,1193,2775,2768)$ for an odd edge.

(viii) We obtain
$935=5\times
11\times 17=468^2-467^2=96^2-91^2=48^2-37^2=36^2-19^2=5(94^2-93^2)=5(14^2-3^2)=11(43^2-42^2)=11(11^2-6^2)=17(28^2-27^2)=17(8^2-3^2)=55(9^2-8^2)=85(6^2-5^2)=187(3^2-2^2)$. But
only $935=96^2-91^2=17(28^2-27^2)$ gives
$96^291^2+17^228^227^2=15540^2=d^2$. Then, we have the next smallest Euler
bricks second type as $(935,17472,25704,17497,25721,31080)$ for an odd edge.

\subsection{Pythagorean Conjecture 9: Euler
bricks third type}
\label{pythagorean_conjecture9}
Explore all odd $n\in \mathbb{N}$ which satisfies
\begin{eqnarray}
  n=a(e^2-f^2)=b(g^2-h^2)~\text{and}~a^2e^2f^2+b^2g^2h^2=d^2
  \label{pythagorean_conjecture9.1}
\end{eqnarray}
where $a,b,d,e, f, g, h~(>1)\in \mathbb{N}$. You can also look for
maximum number of occurrences
of $$[\{(a,e,f),(b,g,h)\},~\{aef,bgh,d\}]$$ for any particular $n$ and
hence maximum occurrence for any $n$. Thus for every such $n$,
satisfying equation (\ref{pythagorean_conjecture7.1}), we can discover
one primitive Euler brick of the third type as
$\{a(e^2-f^2),2aef,2bgh,a(e^2+f^2),b(g^2+h^2),2d\}$ having
$(a(e^2-f^2),2aef,2bgh)$ as its edges; $(a(e^2+f^2),b(g^2+h^2),2d)$ as
its face diagonals \cite{Maiti2024}.

\subsubsection{Example:}
(i) We have
$187=11\times
17=94^2-93^2=14^2-3^2=11(9^2-8^2)=17(6^2-5^2)$. But
only $187=11(9^2-8^2)=17(6^2-5^2)$ gives
$11^29^28^2+17^26^25^2=942^2=d^2$. Thus we have the smallest Euler brick
third type $(187,1584,1020,1595,1037,1884)$ for odd edge.

(ii) We have
$429=3\times 11\times
13=215^2-214^2=73^2-70^2=25^2-14^2=23^2-10^2=3(72^2-71^2)=3(12^2-1^2)=11(20^2-19^2)=11(8^2-5^2)=13(17^2-16^2)=13(7^2-4^2)=33(7^2-6^2)=39(6^2-5^2)=143(2^2-1^2)$. But
only $429=11(8^2-5^2)=39(6^2-5^2)$ gives
$11^28^25^2+39^26^25^2=1250^2=d^2$. Thus we have the next smallest Euler brick
third type $(429,880,2340,979,2379,2550)$ for odd edge.

(iii) We have
$693=3^2\times 7\times
11=347^2-346^2=53^2-46^2=43^2-34^2=37^2-26^2=3(116^2-115^2)=3(40^2-47^2)=3(20^2-13^2)=3(16^2-5^2)=7(50^2-49^2)=7(10^2-1^2)=9(39^2-38^2)=9(9^2-2^2)=11(32^2-31^2)=11(8^2-1^2)=21(17^2-16^2)=21(7^2-4^2)=33(11^2-10^2)=33(5^2-2^2)=63(6^2-5^2)=77(5^2-4^2)=99(4^2-3^2)$. But
only $693=3(16^2-5^2)=7(10^2-1^2)$ gives
$3^216^25^2+7^210^21^2=250^2=d^2$. Thus we have the next smallest Euler brick
third type $(693,480,140,843,707,500)$ for odd edge.

\section{Biquadratic Diophantine Equation Conjecture for perfect cuboids}
\label{biquadratic_conjecture}
\subsection{Biquadratic Conjecture 1: perfect cuboid first and second type}
\label{biquadratic_conjecture1}
Provide all solutions of the biquadratic Diophantine equation
\begin{eqnarray}
  P^4+Q^4+R^4+S^4=T^2  \label{biquadratic_conjecture1.1}
\end{eqnarray}
where odd $P,Q,R,S~(>1)\in \mathbb{N}$ with $gcd(P,Q)=1=gcd(R,S)$ and
even $T\in \mathbb{N}$ including $T=U^2+V^2$ with $gcd(U,V)=1$ and
$T=d(U^2+V^2)$ with $gcd(U,V)=1$ where odd $d,U,V~(>1)\in
\mathbb{N}$. If further $T=U^2+V^2~\&~PQ=RS=UV$ holds, we can discover
a perfect cuboid first type and for $T=d(U^2+V^2)~\&~PQ=RS=dUV$, we
can generate a perfect cuboid second type \cite{Maiti2024}.

\subsection{Biquadratic Conjecture 2: perfect cuboid third and fifth type}
\label{biquadratic_conjecture2}
Find all solutions of the biquadratic Diophantine equation
\begin{eqnarray}
  P^4+Q^4+a^2(R^4+S^4)=T^2  \label{biquadratic_conjecture2.1}
\end{eqnarray}
where odd $a,P,Q,R,S~(>1)\in \mathbb{N}$ with $gcd(P,Q)=1=gcd(R,S)$
and even $T\in \mathbb{N}$ including $T=U^2+V^2$ with $gcd(U,V)=1$ and
$T=d(U^2+V^2)$ with $gcd(U,V)=1$ where odd $d,U,V~(>1)\in
\mathbb{N}$. If further $T=U^2+V^2~\&~PQ=aRS=UV$ holds, we can detect
a perfect cuboid third type and for $T=d(U^2+V^2)~\&~PQ=aRS=dUV$, we
can produce a perfect cuboid fifth type \cite{Maiti2024}.

\subsection{Biquadratic Conjecture 3: perfect cuboid fourth and sixth type}
\label{biquadratic_conjecture3}
Discover all solutions of the biquadratic Diophantine equation
\begin{eqnarray}
  a^2(P^4+Q^4)+b^2(R^4+S^4)=T^2  \label{biquadratic_conjecture3.1}
\end{eqnarray}
where odd $a,b,P,Q,R,S~(>1)\in \mathbb{N}$ with $gcd(P,Q)=1=gcd(R,S)$
and even $T\in \mathbb{N}$ including $T=U^2+V^2$ with $gcd(U,V)=1$ and
$T=d(U^2+V^2)$ with $gcd(U,V)=1$ where odd $d,U,V~(>1)\in
\mathbb{N}$. Further, when $T=U^2+V^2~\&~aPQ=bRS=UV$ satisfies, a
perfect cuboid of the fourth type can be produced and for
$T=d(U^2+V^2)~\&~aPQ=bRS=dUV$, a perfect cuboid of sixth type can be
discovered \cite{Maiti2024}.

\subsection{Remark:}
We can develop similar type of conjectures for even edge.

The research topic on perfect cuboid and Euler brick appeared several
times in American Mathematical Monthly articles and was even a
research topic for a PhD theses in 2000 in Europe and 2004 in
China. This research topic has certainly great educational value
because of its simplicity. As an example, This research topic appeared
in a journal run by undergraduates similar to HCMR \cite{Ortan}.

Noam Elkies told Oliver Knill \cite{Knill}: ``The algebraic surface
parametrizing Euler bricks is the intersection in $P^5$ of the quadrics
$x^2+y^2=c^2$, $z^2+x^2=b^2$, $y^2+z^2=a^2$ which happens to be a $K3$
surface of maximal rank, so quite closely related to much of Knill's
own recent work in number theory. Adding the condition
$x^2+y^2+z^2=d^2$ yields a surface of general type, so it might well
have no nontrivial rational points but nobody knows how to prove such
a thing.'' Noam also remarked that \cite{Knill}: ``Eulers
parametrization would only lead to a finite number of perfect Cuboids
as a consequence of Mordells theorem. There seems however no reason to
be known which would tell whether there are maximally finitely many
primitive perfect cuboids''.

\section{Conclusion}
All the perfect cuboids (if any) belong to the the solutions of first
six conjectures (and hence consequently to the the solutions of last
three conjectures) and all the Euler bricks can be generated from the
the solutions of the last three conjectures \cite{Maiti2024}. Also,
perfect cuboid (if any) can be discovered from these three types of
Euler bricks only. Moreover, all the Euler bricks (body cuboid);
available in Table given by Maurice Kraitchik, Rathbun and Granlund as
well as Rathbun\cite{Rathbun1,Rathbun2,Rathbun3}; can be classified as
the Euler bricks first type, Euler bricks second type and Euler bricks
third type. In the range of odd edge less than 1000, there are 11
primitive Euler bricks. These are 8 Euler bricks second type, 3 Euler
bricks third type and no Euler bricks of first type.

\vspace{1cm}

{\bf Acknowledgment:} {\it I am grateful to the University Grants
  Commission (UGC), New Delhi for awarding the Dr. D. S. Kothari Post
  Doctoral Fellowship from 9th July, 2012 to 8th July, 2015 at Indian
  Institute of Technology (BHU), Varanasi. The self training for this
  investigation was continued during the period. I am also
  acknowledge gratitude to my wife Dr. Ankita Chaturvedi who inspired
  me to complete this study.}

\end{document}